\documentstyle[twoside,psfig]{article}

\oddsidemargin=\evensidemargin 
\catcode`\@=11 \long\def\@makefntext#1{ \protect\noindent \hbox to
3.2pt {\hskip-.9pt
$^{{\eightrm\@thefnmark}}$\hfil}#1\hfill}       

\def\ps@myheadings{\let\@mkboth\@gobbletwo      
\def\@oddhead{\hbox{}
\rightmark\hfil\eightrm\thepage}
\def\@oddfoot{}\def\@evenhead{\eightrm\thepage\hfil
\leftmark\hbox{}}\def\@evenfoot{}
\def\sectionmark##1{}\def\subsectionmark##1{}}

\catcode`@=11                        
\def\ps@plain{\let\@mkboth\@gobbletwo
     \def\@oddhead{}\def\@oddfoot{\eightrm\hfil\thepage
     \hfil}\def\@evenhead{}\let\@evenfoot\@oddfoot}

\newcounter{sectionc}\newcounter{subsectionc}\newcounter{subsubsectionc}
\renewcommand{\section}[1] {\vspace{12pt}\addtocounter{sectionc}{1}
\setcounter{subsectionc}{0}\setcounter{subsubsectionc}{0}\noindent
    {\tenbf\thesectionc. #1}\par\vspace{5pt}}
\renewcommand{\subsection}[1] {\vspace{12pt}\addtocounter{subsectionc}{1}
    \setcounter{subsubsectionc}{0}\noindent
    {\bf\thesectionc.\thesubsectionc.
    {\kern1pt \bfit #1}}\par\vspace{5pt}}
\renewcommand{\subsubsection}[1] {\vspace{12pt}
    \addtocounter{subsubsectionc}{1}
    \noindent
    {\tenrm\thesectionc.\thesubsectionc.\thesubsubsectionc. {\kern1pt
    \it #1}}\par\vspace{5pt}}

\newcounter{appendixc}
\newcounter{subappendixc}[appendixc]
\newcounter{subsubappendixc}[subappendixc]

\renewcommand{\appendix}[1] {\vspace{12pt}  
    \refstepcounter{appendixc}      
    \setcounter{figure}{0}
    \setcounter{table}{0}
    \setcounter{lemma}{0}
    \setcounter{theorem}{0}
    \setcounter{corollary}{0}
    \setcounter{definition}{0}
    \setcounter{equation}{0}
    \renewcommand{\thefigure}{\Alph{appendixc}.\arabic{figure}}
    \renewcommand{\thetable}{\Alph{appendixc}.\arabic{table}}
    \renewcommand{\theappendixc}{\Alph{appendixc}}
    \renewcommand{\thelemma}{\Alph{appendixc}.\arabic{lemma}}
    \renewcommand{\thetheorem}{\Alph{appendixc}.\arabic{theorem}}
    \renewcommand{\thedefinition}{\Alph{appendixc}.\arabic{definition}}
    \renewcommand{\thecorollary}{\Alph{appendixc}.\arabic{corollary}}
    \renewcommand{\theequation}{\Alph{appendixc}.\arabic{equation}}
    \noindent{\tenbf Appendix \theappendixc #1}\par\vspace{5pt}}

\newcommand{\smalllineskip}{\baselineskip=10pt}

\newcommand{\copyrightheading}[1]
    {\vspace*{-2.5cm}\smalllineskip{\flushleft
    {\footnotesize }\\
    {\footnotesize \copyright\kern2pt }\\
         }}

\newcounter{itemlistc}
\newcounter{romanlistc}
\newcounter{alphlistc}
\newcounter{arabiclistc}

\newcommand{\fcaption}[1]{
        \refstepcounter{figure}
        \setbox\@tempboxa = \hbox{\footnotesize Fig.~\thefigure. #1}
        \ifdim \wd\@tempboxa > 5in
           {\begin{center}
        \parbox{5in}{\footnotesize\smalllineskip Fig.~\thefigure. #1}
            \end{center}}
        \else
             {\begin{center}
             {\footnotesize Fig.~\thefigure. #1}
              \end{center}}
        \fi}

\newcommand{\tcaption}[1]{
        \refstepcounter{table}
        \setbox\@tempboxa = \hbox{\footnotesize Table~\thetable. #1}
        \ifdim \wd\@tempboxa > 5in
           {\begin{center}
        \parbox{5in}{\footnotesize\smalllineskip Table~\thetable. #1}
            \end{center}}
        \else
             {\begin{center}
             {\footnotesize Table~\thetable. #1}
              \end{center}}
        \fi}

\def\pmb#1{\setbox0=\hbox{#1}
    \kern-.025em\copy0\kern-\wd0
    \kern.05em\copy0\kern-\wd0
    \kern-.025em\raise.0433em\box0}

\def\fnt#1#2{\footnotetext{\kern-.3em
    {$^{\mbox{\scriptsize #1}}$}{#2}}}

\def\fpage#1{\begingroup
\voffset=.3in
\thispagestyle{empty}\begin{table}[b]\centerline{\footnotesize #1}
    \end{table}\endgroup}

\def\runninghead#1#2{\pagestyle{myheadings}
\markboth{{\protect\footnotesize\it{\quad #1}}\hfill}
{\hfill{\protect\footnotesize\it{#2\quad}}}}

\font\tenrm=cmr10  \font\tenbf=cmbx10
\font\bfit=cmbxti10 at 10pt \font\ninerm=cmr9 
 \font\eightrm=cmr8

\newtheorem{theorem}{Theorem}   

\newtheorem{definition}{Definition}

\def\@begintheorem#1#2{\trivlist    
    \item[\hskip\labelsep{\bf #1\ #2.}]}
\def\@opargbegintheorem#1#2#3{\trivlist
    \item[\hskip\labelsep{\bf #1\ #2\ (#3).}]}

        {\setcounter{itemlistc}{0}      
     \begin{list}{$\bullet$}        
    {\usecounter{itemlistc}         
     \leftmargin10pt           
     \setlength{\parsep}{0pt}
     \setlength{\itemsep}{0pt}     
    }}{\end{list}}

    {\setcounter{romanlistc}{0}     
     \begin{list}{$($\roman{romanlistc}$)$} 
    {\usecounter{romanlistc}        
     \leftmargin18pt 
     \setlength{\parsep}{0pt}
     \setlength{\itemsep}{0pt}  
     \settowidth{\labelwidth}{#1}
    }}{\end{list}}

    {\setcounter{enumii}{0}         
     \begin{list}{$($\alph{enumii}$)$}  
    {\usecounter{enumii}            
     \leftmargin18pt        
     \setlength{\parsep}{0pt}
     \setlength{\itemsep}{0pt}  
     \settowidth{\labelwidth}{#1}
    }}{\end{list}}

\def\qed{\hbox{${\vcenter{\vbox{            
   \hrule height 0.4pt\hbox{\vrule width 0.4pt height 6pt
   \kern5pt\vrule width 0.4pt}\hrule height 0.4pt}}}$}}

\def\theequation{\thesectionc.\arabic{equation}}  
\begin{document}

\runninghead{Hyperbolic Volume of Link Families} {Hyperbolic Volume
of Link Families}

\setcounter{page}{1}

\markboth{Slavik Jablan and Ljiljana Radovi\' c}{}


\fpage{1} \centerline{\bf HYPERBOLIC VOLUME OF LINK FAMILIES}
\bigskip

\centerline{\footnotesize SLAVIK JABLAN and LJILJANA RADOVI\'
C$^\dag $}
\medskip
\centerline{\footnotesize\it The Mathematical Institute, Knez
Mihailova 35,}\centerline{\footnotesize\it P.O.Box 367, 11001
Belgrade,}\centerline{\footnotesize\it Serbia}
\centerline{\footnotesize\it jablans@mi.sanu.ac.yu}

\medskip

\centerline{\footnotesize\it Faculty of Mechanical Engineering$^\dag
$}\centerline{\footnotesize\it A.~Medvedeva 14 }
\centerline{\footnotesize\it 18 000 Ni\v s }
\centerline{\footnotesize\it Serbia} \centerline{\footnotesize\it
liki@masfak.ni.ac.yu}

\bigskip

\begin{abstract}
For families of knots and links given in Conway notation we compute
lower maximal and upper minimal bound of hyperbolic volume by using
source links and augmented links.
\end{abstract}

In this paper we consider families of knots and links ($KL$s) given
in Conway notation and their hyperbolic volume. Hyperbolic volume of
the complement of a link $L$ will be shortly denoted as $Vol(L)$.
First we define a $KL$ family [1]:

\begin{definition}
For a link or knot $L$ given in an unreduced\footnote{The Conway
notation is called unreduced if in symbols of polyhedral $KL$s
elementary tangles 1 in single vertices are not omitted.} Conway
notation $C(L)$ denote by $S$ a set of numbers in the Conway symbol
excluding numbers denoting basic polyhedra and zeros (denoting the
position of tangles in the vertices of polyhedra). For $C(L)$ and an
arbitrary (non-empty) subset $\tilde S$ of $S$ the family $F_{\tilde
S}(L)$ of knots or links derived from $L$ is constructed by
substituting each $a \in S_f$, $a\neq 1$, by $sgn(a) (|a|+k_a)$ for
$k_a \in N$.
\end{definition}

\begin{definition}
A $KL$ with single bigons, or equivalently, a $KL$ given by Conway
symbol containing only tangles $1$, $-1$, $2$, or $-2$ is called a
{\it source link}.
\end{definition}

This means that all $KL$s generated from a source link $S$ by
substituting single bigons by chains of bigons make a {\it family}
generated from $S$.

\begin{definition}
Let $D$ be a reduced link diagram. Two crossings of $D$ are twist
equivalent if there exist a flype connecting these two crossings
into a bigon. The {\it twist number of diagram} $D$ is the number of
its twist equivalence classes. The {\it twist number}
$t(L)=t_{min}(L)$ of a link $L$ is the minimal twist number over all
diagrams of $L$ [2,3].
\end{definition}

Since Conway symbols of links are twist-reduced, for a link diagram
given in Conway notation twist number $t_D$ is the number of
parameters plus the number of single (isolated) crossings in the
Conway symbol. All diagrams in a family of alternating link diagrams
have the same twist number $t_D$. For example, twist number $t_D$ of
all diagrams of the family $p\,q$ ($p\ge 2$, $q\ge 2$) is
$t_D(p\,q)=2$, for the family of the diagrams $p\,1\,q$ ($p\ge 2$,
$q\ge 2$) it is $t_D(p\,1\,q)=3$, {\it etc.} However, there are
alternating links with non-alternating diagrams of smaller twist
number than the twist number of their alternating diagrams [2].

\bigskip

\noindent {\bf Conjecture 1.} For every link $L$ given by a minimal
twist-reduced diagram $D$, the twist number $t(L)\ge {\lfloor {t_D
\over 2} \rfloor}+1$.

\bigskip

In many papers the twist number $t(L)$ is used as the main tool to
determine the upper and lower bound of a hyperbolic volume.

In the book {\it LinKnot} [3] it is pointed out  that different knot
and link ($KL$) invariants (unlinking number, signature,
coefficients of Alexander and Jones polynomial, {\it etc.}) can be
expressed as the functions of parameters from Conway symbols of $KL$
families. Experimental results suggest that the same holds for
hyperbolic volume.

Hyperbolic volumes of families of $KL$s  from the program {\it
LinKnot}, given in Conway notation for $KL$s up to 49 crossings, are
computed by {\it Knotscape} and Windows version of {\it Snap Pea},
so all computations are done with a limited precision.

Simplest one-parameter family of knots is the family $2p+1$ ($p\ge
1$) which consists of non-hyperbolic knots $3_1$, $5_1$, $\ldots $.

\begin{figure}[th]
\centerline{\psfig{file=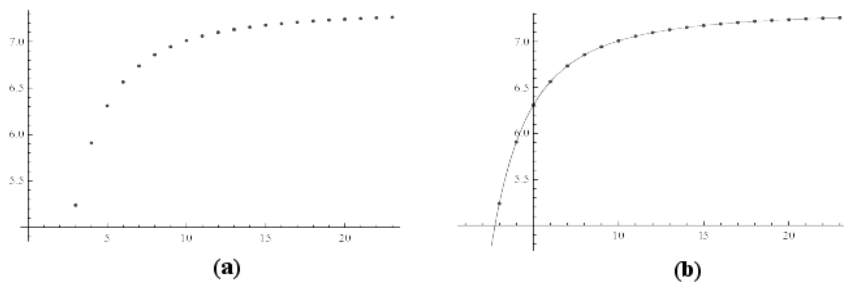,width=4.60in}} \vspace*{8pt}
\caption{(a) Point plot of the hyperbolic volumes for the family
$(p+1)\,(p+1)$ ($p\ge 1$); (b) interpolation  by a function given by
equation 1.1 for $n=4$.\label{f1.1}}
\end{figure}

For $KL$s from one-parameter subfamily $(p+1)\,(p+1)$ ($p\ge 1$) of
two-parameter link family  $p\,q$ ($p\ge 2$, $q\ge 2$) with $p\le
23$, hyperbolic volumes are given in the following table:

\bigskip

\begin{tabular}{|c|c|c|c|} \hline

$2\,2$ &  $2.0298832128$ & $14\,14$ & $7.1290060758$ \\  \hline

$3\,3$ &  $4.059766426$ & $15\,15$ & $7.154316936$ \\  \hline

$4\,4$ &  $5.2386841008$ & $16\,16$ & $7.1750981657$ \\  \hline

$5\,5$ &  $5.907963404$ &  $17\,17$ & $7.192366348$ \\  \hline

$6\,6$ &  $6.3090903924$ &  $18\,18$ & $7.2068688595$ \\  \hline

$7\,7$ &  $6.564505128$ &  $19\,19$ & $7.219164881$ \\  \hline

$8\,8$ &  $6.7359047525$ & $20\,20$ & $7.2296794015$ \\  \hline

$9\,9$ &  $6.856023126$ & $21\,21$ & $7.238740004$ \\  \hline

$10\,10$ &  $6.9432605638$ &  $22\,22$ & $7.2466024333$ \\  \hline

$11\,11$ &  $7.008519846$ &  $23\,23$ & $7.253468667$ \\  \hline

$12\,12$ &  $7.0585637385$ & $24\,24$ & $7.2594999144$ \\  \hline

$13\,13$ & $7.097755265$ & $$ & $$  \\  \hline

\end{tabular}

\bigskip
The list plot (Fig. 1a) was obtained by {\it Mathematica 6.0}. All
discrete list plots for one-parameter families are interpolated by
functions of the form
$${{\sum_{i=0}^{n}a_ix^{2i}} \over {\sum_{i=0}^{n}b_ix^{2i}}}+c
\eqno{(1.1)}$$ with $n=3$ or $n=4$. In this example, the best
interpolation is obtained by {\it Mathematica} function {\bf
FindFit} for $n=4$, with the maximal error $2.17782\times 10^{-9}$
(Fig. 1b). The other possible interpolation by simpler functions of
the form

$$a_0+{\sum_{i=1}^{n}{a_i \over (x+c)^{2i}}} \eqno{(1.2)}$$

\noindent in some cases is less precise.

From the functions interpolating hyperbolic volumes of families we
can make conclusion about their asymptotic behavior. For example, in
the interpolating function for the family $p\,p$ ($p\ge 2$) obtained
for $n=4$ the coefficients are $a_8=2.3491324728718244$,
$b_8=0.5358879857172603$, and $c=2.944097878883564$, so the
interpolating function converges to ${a_8 \over b_8} +c
=Vol(6^*)=7.32772...$, which is the hyperbolic volume of Borromean
rings complement. This result is in the complete agreement with the
results of C.~Petronio and A.~Vesnin [4, Proposition 1], which
showed that hyperbolic volume of the family $p\,q$ ($p\ge 2$, $q\ge
2$) converges to $V_2=Vol(6^*)=7.32772...$ when $p\rightarrow \infty
$ and $q\rightarrow \infty $. The same result we obtained for the
family $p\,1\,p$ ($p\ge 2$), i.e., for the family $p\,1\,q$ when
$p\rightarrow \infty $ and $q\rightarrow \infty $.

For the families of rational $KL$s of the form $p\,\ldots \,p$
($p\ge 2$) where $p$ occurs $n$ times, shortly denoted as $p^n$,
$lim_{p\rightarrow \infty} Vol(p^n)=(n-1)V_2$. The same holds for
the rational link families of the form $p\,1\,p$, $p\,1\,p\,1\,p$,
$p\,1\,p\,1\,p\,1\,p$, ... where $p$ occurs $n$ times.

Augmented links are used for obtaining bounds of hyperbolic volumes
[5,6,7]. In the language of Conway symbols and chains of bigons,
augmentation of a bigon chain (tangle) $p$ ($p\ge 2$) is its
replacement by tangle $(2,-2)\,0$ (Fig. 2).

\begin{figure}[th]
\centerline{\psfig{file=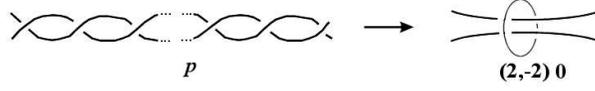,width=3.20in}} \vspace*{8pt}
\caption{Augmentation of a tangle $p$ ($p\ge 2$). \label{f2.1}}
\end{figure}

\bigskip

\begin{definition}
For a family of alternating links $L$ given by the Conway symbol
$$C(p_1,...,p_k,...p_n)$$ containing chain of bigons $p_k$ ($p_k\ge
2$), subfamily with changing parameter $p_k$ ($p_k=2,3,4,\ldots$)
and all other $p$-s fixed is called {\it $p_k$-subfamily}, link
obtained by replacing $p_k$ ($k\in \{1,2,\ldots,n\}$, $p_k\ge 2$) by
a tangle $(2,-2)\,0$ is called {\it $p_k$-augmented link} and
denoted by $C(p_1,...,\overline p_k,...p_n)$, and link obtained by
replacing $p_k$ ($k\in \{1,2,\ldots,n\}$, $p_k\ge 2$) by a single
bigon $2$ is called {\it $p_k$-source link} and denoted by
$C(p_1,...,\underline p_k,...p_n)$. In the replacements single
crossings remain unchanged.
\end{definition}

\bigskip

\begin{theorem}
For every alternating algebraic link $L$ with at least two chains of
bigons, or polyhedral link with at least one chain of bigons,
$$Vol(C(p_1,...,\underline
p_k,...p_n))\le Vol(C(p_1,...,p_k,...p_n))<$$
$$Vol(C(p_1,...,p_k+1,...p_n))\le Vol(C(p_1,...,\overline
p_k,...p_n)).$$
\end{theorem}

\bigskip

In general, $$Vol(\underline L) \le Vol(L) \le
 Vol(\overline L),$$
where $\underline L$ is the source link, $\overline L$ is completely
augmented link with all chains of bigons replaced by tangles
$(2,-2)\,0$ and
$$lim Vol(C(p_1,...,p_n))=\overline L,$$ where the limes is taken
when all chains of bigons tend to infinity [5, Corollary 2].

\bigskip

For example, let's consider two-parameter family $L=8^*p\,0.q\,0$
($p\ge 2$, $q\ge 2$). For fixed $q=2,3,4,5\ldots $ we obtain
sequence of one-parameter subfamilies: $f_2=8^*p\,0.2\,0$,
$f_3=8^*p\,0.3\,0$, $f_4=8^*p\,0.4\,0$, $f_5=8^*p\,0.5\,0$, ...
satisfying the relations:

$$16.6380380564\ldots =Vol(8^*2\,0.2\,0)\le Vol(f_2) \le
Vol(8^*(2,-2).2\,0)=19.29865114\ldots $$

$$17.7392681473\ldots =Vol(8^*2\,0.3\,0)\le Vol(f_3) \le
Vol(8^*(2,-2).3\,0)=20.559914\ldots $$

$$18.3010568281\ldots =Vol(8^*2\,0.4\,0)\le Vol(f_4) \le
Vol(8^*(2,-2).4\,0)=21.21212466\ldots $$

$$18.6120521177\ldots =Vol(8^*2\,0.5\,0)\le Vol(f_5) \le
Vol(8^*(2,-2).5\,0)=21.5747527\ldots $$

\noindent From the interpolating functions for $f_2$, $f_3$, $f_4$,
$f_5$ we obtain asymptotic values: $19.2972\ldots $, $20.5586\ldots
$, $21.2112\ldots $, $21.574\ldots $, respectively.

For the whole family $8^*p\,0.q\,0$, from the subfamily
$8^*p\,0.p\,0$ and interpolating function 1.1 with $n=4$,
$a_8=17.378499561645388$, $b_8=10.894820228608358$,
$c=20.771542391115194$, we obtain the asymptotic value ${a_8 \over
b_8} +c =22.3667\ldots $, and $Vol(\overline
L)=Vol(8^*(2,-2).(2,-2))=22.36710548...$ Hence, for the family
$L=8^*p\,0.q\,0$ ($p\ge 2$, $q\ge 2$) holds the following relation:

$$16.6380380564\ldots =Vol(\underline L)\le Vol{L}\le Vol(\overline
L)=22.36710548\ldots $$

Conjecture 2 can be applied ed also to non-alternating $KL$s, but
only positive chains of bigons can be varied and augmented. For
example, let's consider two-parameter family of non-alternating
knots given by minimal diagrams $L=8^*p\,0.-q\,0$ ($p\ge 2$, $q\ge
2$). For fixed $q=2,3,4,5\ldots $ we obtain sequence of
one-parameter subfamilies: $f_2'=8^*p\,0.-2\,0$,
$f_3'=8^*p\,0.-3\,0$, $f_4'=8^*p\,0.-4\,0$, $f_5'=8^*p\,0.-5\,0$,
... satisfying the relations:

$$13.2900030686\ldots =Vol(8^*2\,0.-2\,0)\le Vol(f_2') \le
Vol(8^*(2,-2).-2\,0)=16.69568447\ldots $$

$$17.7392681473\ldots =Vol(8^*2\,0.-3\,0)\le Vol(f_3') \le
Vol(8^*(2,-2).-3\,0)=19.29865114\ldots $$

$$18.3010568281\ldots =Vol(8^*2\,0.-4\,0)\le Vol(f_4') \le
Vol(8^*(2,-2).-4\,0)=20.559914\ldots $$

$$18.6120521177\ldots =Vol(8^*2\,0.-5\,0)\le Vol(f_5') \le
Vol(8^*(2,-2).-5\,0)=21.21212466\ldots $$

\noindent From the interpolating functions for $f_2$, $f_3$, $f_4$,
$f_5$ we obtain asymptotic values: $16.6957\ldots $, $19.2961\ldots
$, $20.56\ldots $, $21.2093\ldots $, respectively.



However, we cannot obtain correct results by fixing $p$, and varying
$q$. This means that we need to fix negative chains and vary only
positive ones. For negative chains, this problem can be solved by
taking mirror image of a link $L$ or another minimal diagram of the
same link with positive chains of bigons..

For completely augmented links signs of bigon chains are not
relevant. For example, the family $L=8^*p\,0.-q\,0$ ($p\ge 2$, $q\ge
2$) also can be  given by the minimal diagram
$L'=9^*(p-1)\,0.-1.-1.q\,0.-1.-1:-1.-1$. From completely augmened
link $\overline L'$ we obtain minimal upper bound $\overline L
=\overline L'
=Vol(9^*(2,-2).-1.-1.(2,-2).-1.-1:-1.-1)=22.36710548\ldots $, so we
conclude that the family $L=8^*p\,0.-q\,0$ ($p\ge 2$, $q\ge 2$)
holds the following relation:

$$13.2900030686\ldots =Vol(\underline L)\le Vol{L}\le Vol(\overline
L)=22.36710548\ldots $$

\noindent Hence, the families $8^*p\,0.q\,0$ and $8^*p\,0.-q\,0$
($p\ge 2$, $q\ge 2$) have the same minimal upper bound.

In the same way, we can consider the minimal diagrams
$8^*p.-1.q.-1.r.-1:-1$ and $.q.-p.-(r+1).2\,0:-1$ ($p,q,r\ge 2$) of
the same non-alternating link. Their corresponding completely
augmented links $8^*(2,-2)\,0.-1.(2,-2)\,0.-1.(2,-2)\,0.-1:-1$ and
$.(2,-2)\,0.(2,-2)\,0.(2,-2)\,0.2\,0:-1$ have the same hyperbolic
volume $18.83168337\ldots$ (Fig. 3).

\begin{figure}[th]
\centerline{\psfig{file=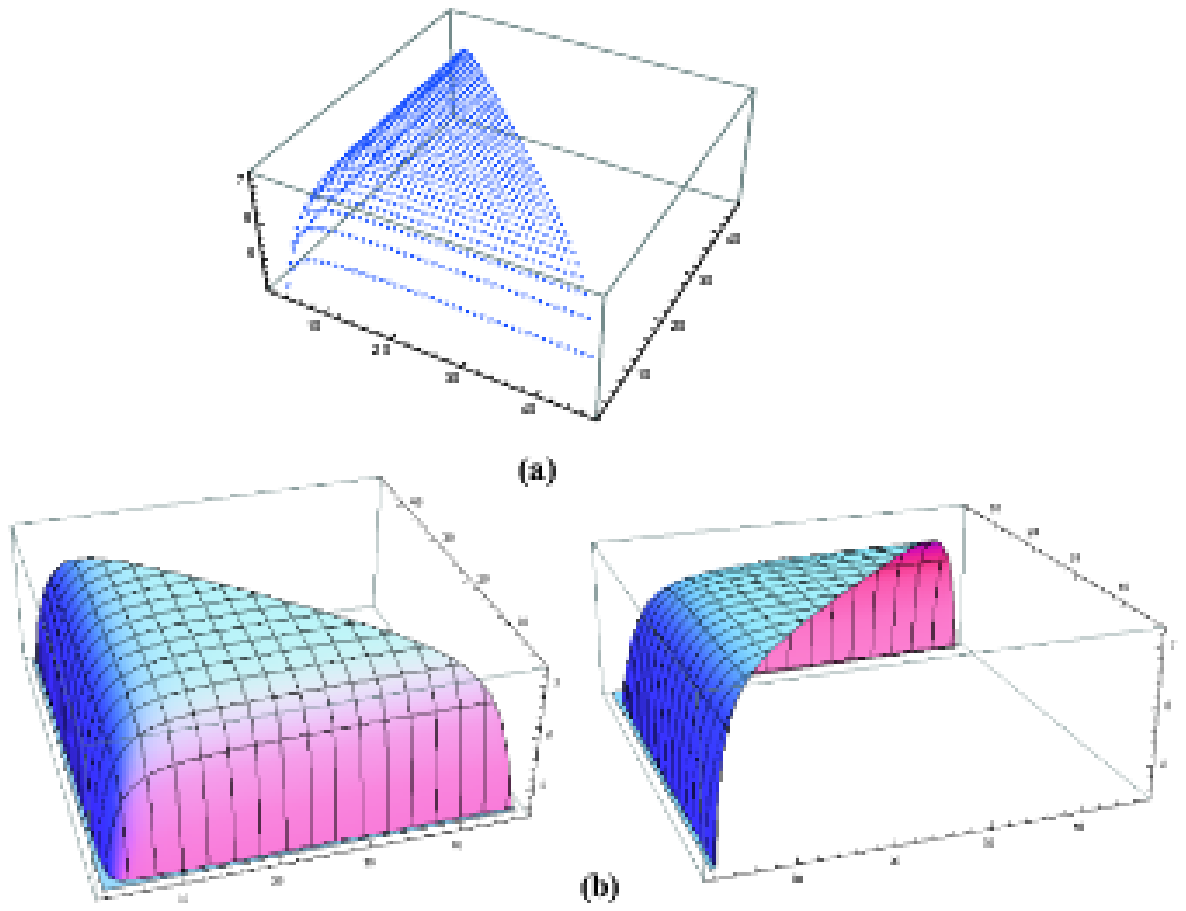,width=3.80in}} \vspace*{8pt}
\caption{(a) Point plot  and (b) list plot 3D of the hyperbolic
volumes for the family $p\,q$ ($p\ge 2$), $q\ge 2$. \label{f3.1}}
\end{figure}

\begin{figure}[th]
\centerline{\psfig{file=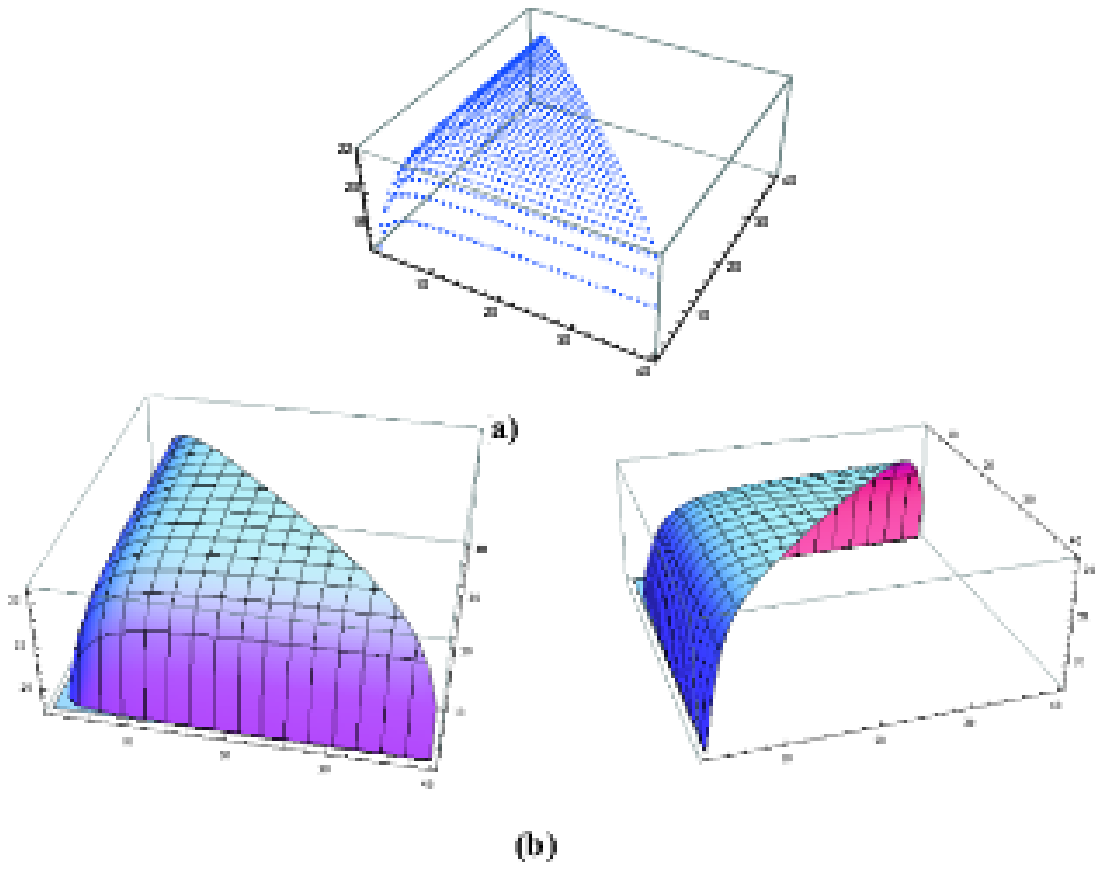,width=3.70in}} \vspace*{4pt}
\caption{(a) Point plot  and (b) list plot 3D of the hyperbolic
volumes for the family $8^*p\,0.q\,0$ ($p\ge 2$, $q\ge 2$).
\label{f4.1}}
\end{figure}

\newpage

As the last example of this kind we can consider three families of
knots given by minimal diagrams $10^*p\,0::.q\,0$,
$10^*-p\,0::.q\,0$, $10^*-p\,0::.-q\,0$ ($p\ge 2$, $q\ge 2$).
Despite of the fact that the first is the family of alternating
links, and other two are non-alternating, their hyperbolic volumes
converge to the hyperbolic volume of the same completely augmented
knot $10^*(2,-2)::.(2,-2)$ with the hyperbolic volume
$26.3062315\ldots $ (Fig. 5).

\begin{figure}[th]
\centerline{\psfig{file=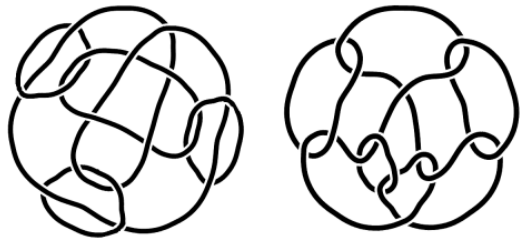,width=2.80in}} \vspace*{8pt}
\caption{Completely augmented links
$8^*(2,-2)\,0.-1.(2,-2)\,0.-1.(2,-2)\,0.-1:-1$ and
$.(2,-2)\,0.(2,-2)\,0.(2,-2)\,0.2\,0:-1$ with the same hyperbolic
volume. \label{f5.1}}
\end{figure}

\begin{figure}[th]
\centerline{\psfig{file=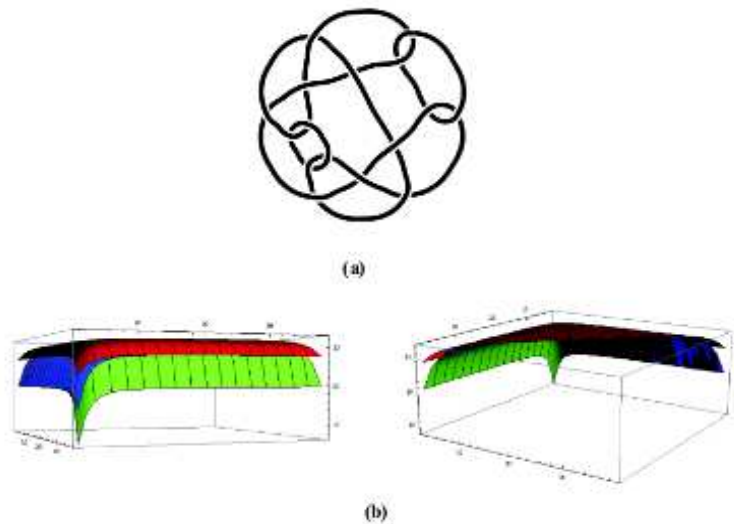,width=4.20in}} \vspace*{8pt}
\caption{(a) Completely augmented knot $10^*(2,-2)::.(2,-2)$; (b)
list plot 3D of the hyperbolic volumes for the families
$10^*p\,0::.q\,0$, $10^*-p\,0::.q\,0$, $10^*-p\,0::.-q\,0$ ($p\ge
2$, $q\ge 2$) converging to the same minimal upper bound.
\label{f6.1}}
\end{figure}

\newpage
Two-parameter families and their hyperbolic volumes are visualized
using {\it Mathematica} functions {\bf ListPointPlot3D} and {\bf
ListPlot3D}. For the $KL$ family $p\,q$ ($p\ge 2$, $q\ge 2$),
beginning with knot $2\,2$ ($4_1$) the point plot ({\bf
ListPointPlot3D}) and the plot ({\bf ListPlot3D}) of the list
$(p,q,Vol(p\,q))$ are shown in Fig. 3. The same data for the family
$8^*p\,0.q\,0$ ($p\ge 2$, $q\ge 2$) are shown on Fig. 4. For all
alternating two-parameter $KL$ families we obtained smooth surfaces
with one-parameter subfamilies which can be interpolated by the
preceding class of functions.

\begin{figure}[th]
\centerline{\psfig{file=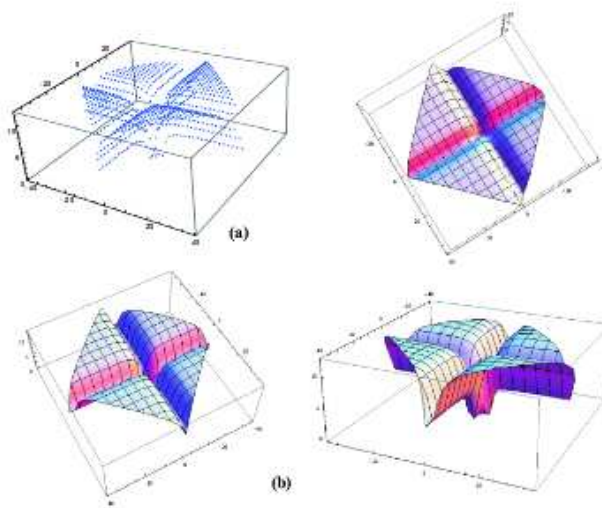,width=3.60in}} \vspace*{8pt}
\caption{(a) Point plot  and (b) list plot 3D of the hyperbolic
volumes for the family $6^*-(2p+1).(2q).-2.2.-2$. \label{f7.1}}
\end{figure}

\begin{figure}[th]
\centerline{\psfig{file=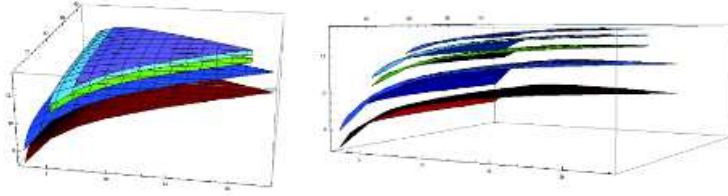,width=4.00in}} \vspace*{8pt}
\caption{List plot 3D of the hyperbolic volumes for the
three-parameter family of pretzel knots $p,q,r$ for fixed $p\in
\{2,3,4,5, 6\}$. \label{f8.1}}
\end{figure}

\begin{figure}[th]
\centerline{\psfig{file=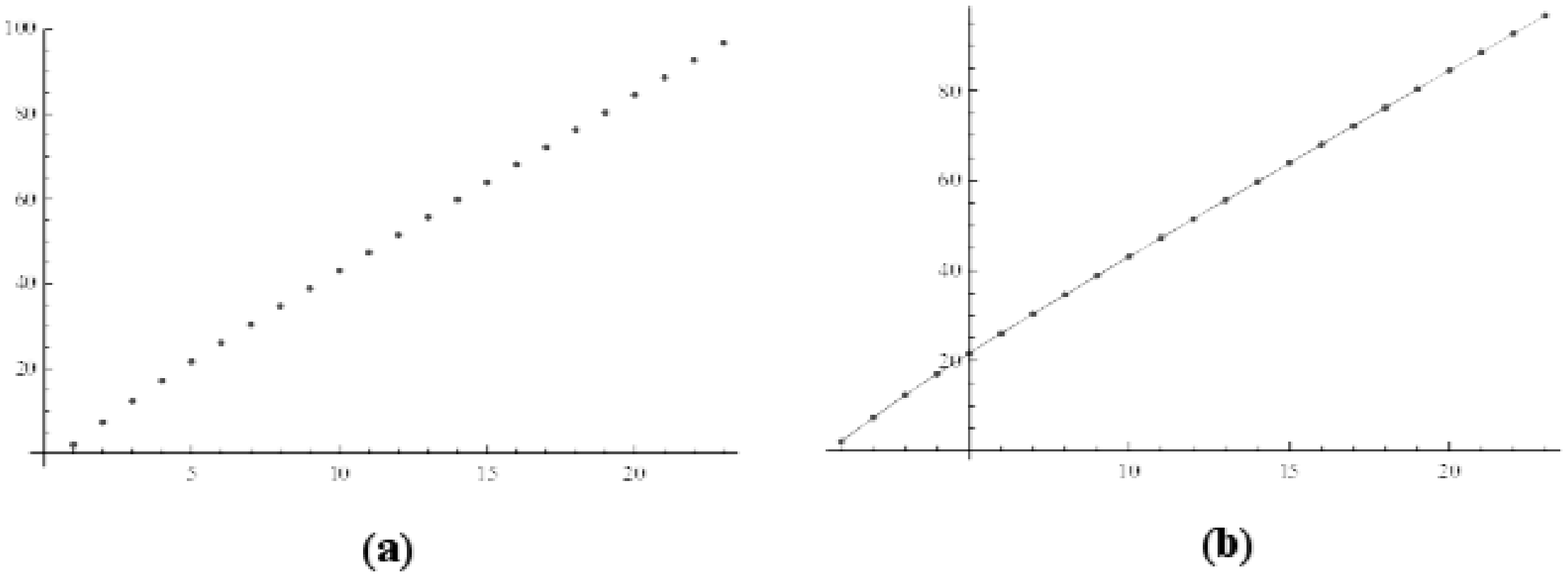,width=3.70in}} \vspace*{8pt}
\caption{(a) Point plot of the hyperbolic volumes for the family of
the antiprismatic basic polyhedra $(2n)^*$ ($n=2,\ldots 24$); (b)
its interpolation function given by equation (1.1)
 for $n=4$.\label{f9.1}}
\end{figure}

Similar results are obtained for two-parameter families of
non-alternating knots, up to 49 crossings. The most interesting of
all are point and list plots of hyperbolic volumes of knots
belonging to the family $6^*-(2p+1).(2q).-2.2.-2$ where both
positive and negative values for $p$ and $q$ are allowed
\footnote{This family contains Lorenz knots for $p\ge 2$, $q\ge
1$, with additional condition $p\le q+1$.}.

In order to visually represent hyperbolic volumes of knots belonging
to three-parameter families, we fix one parameter and vary remaining
two. Fig. 8  shows three-dimensional list plots of hyperbolic
volumes for pretzel knots $p,q,r$ for $p=2$, $p=3$, $p=4$, $p=5$ and
$p=6$ simultaneously, where values of $q$ and $r$ are varied.
Computations include alternating and non-alternating pretzel knots
$p,q,r$ up to 49 crossings.

Similar results are obtained for the hyperbolic volume of the
families of basic polyhedra computed for the family of antiprismatic
basic polyhedra $6^*$, $8^*$,$\ldots $, $48^*$, given by braid words
of the form $(aB)^n$ ($n\ge 2$)\footnote{For $n=2$ we obtain, as a
limiting case, figure-eight knot.} (Fig. 9) and for the family of
basic polyhedra $9^*$ ($AbACbACbC$), $10^{**}$ ($AbAbCbACbC$),
$11^{**}$ ($AbAbACbACbC$), $12F$ ($AbAbAbCbACbC$), {\it etc.}
(computed up to 48 crossings).

The obtained results suggest that hyperbolic volume of $KL$s given
in Conway notation, with arbitrarily large number of crossings can
be computed (or at least approximated) directly from the hyperbolic
volume of the $KL$ family it belongs to provided we have the
interpolating function, and that maximal lower and minimal upper
bound of the hyperbolic volume for a given family can be simply
computed from source link and augmented link.

The lower and upper bounds of hyperbolic volumes for all alternating
$KL$ families derived from source links with at most $n=9$ crossings
are given in the following table. As the main referential volumes
are used volume of ideal hyperbolic tetrahedron
$V_0=1.0149416064...$, hyperbolic volume of Whitehead link
$V_1=Vol(2\,1\,2)=3.663862377...$ and the hyperbolic volume of
Borromean rings $V_2=2V_1=Vol(6^*)=7.327724753...$ For every link
family is given its source link in classical and Conway notation
[8], Conway symbol of the family, and lower and upper bound of
hyperbolic volume. Families which have the same lower and upper
bound are given in pairs or triples. All parameters $p$, $q$,
$\ldots $ are greater or equal 2.

\bigskip

\begin{tabular}{|c|c|c|c|c|c|} \hline

$1)$ & $4_1$ & $2\,2$ & $p\,q$ & $2V_0$ & $2V_1$ \\  \hline

$2)$ & $5_1^2$ & $2\,1\,2$ & $p\,1\,q$ & $V_1$ & $2V_1$ \\  \hline

$3)$ & $6_3^2$ & $2\,2\,2$ & $p\,q\,r$ & $5.3334895669\ldots $ & $4V_1$ \\

$4)$ & $6_1^3$ & $2,2,2$ & $p,q,r$ & $5.3334895669\ldots $ & $4V_1$ \\
\hline

$5)$ & $6_3$ & $2\,1\,1\,2$ & $p\,1\,1\,q$ & $5.6930210913\ldots $ & $4V_1$ \\
\hline

$6)$ & $7_6$ & $2\,2\,1\,2$ & $p\,q\,1\,r$ & $7.0849259535\ldots $ & $4V_1$ \\
\hline

$7)$ & $7_7$ & $2\,1\,1\,1\,2$ & $p\,1\,1\,1\,q$ & $7.643375172\ldots $ & $11.7518362\ldots $ \\
\hline

$8)$ & $7_5^2$ & $2\,1,2,2$ & $p\,1,q,r$ & $7.706911803\ldots $ &
$16.0004687\ldots $
\\

$9)$ & $7_1^3$ & $2,2,2+$ & $p,q,r+$ & $7.706911803\ldots $ &
$16.0004687\ldots $ \\  \hline

$10)$ & $7_6^2$ & $.2$ & $.p$ & $8.997351944\ldots $ & $10.991871\ldots $ \\
\hline

$11)$ & $8_{12}$ & $2\,2\,2\,2$ & $p\,q\,r\,s$ & $8.935856927\ldots $ & $6V_1$ \\
\hline

$12)$ & $8_7^2$ & $2\,1\,2\,1\,2$ & $p\,1\,q\,1\,r$ & $8.830664955\ldots $ & $4V_1$ \\
\hline

$13)$ & $8_8^2$ & $2\,1\,1\,1\,1\,2$ & $p\,1\,1\,1\,1\,q$ & $9.672807731\ldots $ & $13.9396857\ldots $ \\
\hline

$14)$ & $8_1^4$ & $2,2,2,2$ & $p,q,r,s$ & $10.14941606\ldots $ &
$24.09218408\ldots $
\\

$15)$ & $8_4^3$ & $(2,2)\,(2,2)$ & $(p,q)\,(r,s)$ & $10.14941606\ldots $ & $24.09218408\ldots $ \\
\hline

$16)$ & $8_9^2$ & $2\,2,2,2$ & $p\,q,r,s$ & $8.967360849\ldots $ & $6V_1$ \\

$17)$ & $8_3^3$ & $2,2,2++$ & $p,q,r++$ & $8.967360849\ldots $ & $6V_1$ \\
\hline

$18)$ & $8_{10}^2$ & $2\,1\,1,2,2$ & $p\,1\,1\,q,r$ & $9.659498545\ldots $ & $17.6277542\ldots $ \\

$19)$ & $8_{12}^2$ & $2\,1,2,2+$ & $p\,1,q,r+$ & $9.659498545\ldots $ & $17.6277542\ldots $ \\
\hline

$20)$ & $8_{15}$ & $2\,1,2\,1,2$ & $p\,1,q\,1,r$ &  $9.930648294\ldots $ & $17.6277542\ldots $ \\
\hline

$21)$ & $8_{13}^2$ & $.2\,1$ & $.p\,1$ & $11.3707742\ldots $ & $13.81327844\ldots $ \\
\hline

$22)$ & $8_{14}^2$ & $.2:2$ & $.p:q$ & $10.6669791\ldots $ & $4V_1$ \\

$23)$ & $8_6^3$ & $.2:2\,0$ & $.p:q\,0$ & $10.6669791\ldots $ & $4V_1$  \\
\hline

$24)$ & $8_{16}$ & $.2.2\,0$ & $.p.q\,0$ & $10.57902192\ldots $ & $15.03537979\ldots $ \\
\hline

$25)$ & $8_{17}$ & $.2.2$ & $.p.q$ & $10.98590761\ldots $ &  $16.11428997\ldots $ \\
\hline

\end{tabular}

\bigskip

\begin{tabular}{|c|c|c|c|c|c|} \hline

$26)$ & $9_{23}$ & $2\,2\,1\,2\,2$ & $p\,q\,1\,r\,s$ & $10.6113483\ldots $ & $6V_1$ \\
\hline

$27)$ & $9_{11}^2$ & $2\,2\,2\,1\,2$ & $p\,q\,r\,1\,s$ & $10.75904664\ldots $ & $6V_1$ \\
\hline

$28)$ & $9_{27}$ & $2\,1\,2\,1\,1\,2$ & $p\,1\,q\,1\,1\,r$ & $10.99998096\ldots $ & $17.47714082\ldots $ \\
\hline

$29)$ & $9_{12}^2$ & $2\,2\,1\,1\,1\,2$ & $p\,q\,1\,1\,1\,r$ & $11.1884778\ldots $ & $19.0795609\ldots $ \\
\hline

$30)$ & $9_{24}^2$ & $2\,1,2\,1,2\,1$ & $p\,1,q\,1,r\,1$ & $12.046092\ldots $ & $18.8316834\ldots $ \\
\hline

$31)$ & $9_{18}^2$ & $2\,2\,1,2,2$ & $p\,q\,1,r,s$ & $11.3817861\ldots $ & $19.5826692\ldots $ \\
\hline

$32)$ & $9_{25}$ & $2\,2,2\,1,2$ & $p\,q,r\,1,s$ & $11.39030515\ldots $ & $23.3281935\ldots $ \\
\hline

$33)$ & $9_1^3$ & $2\,1\,2,2,2$ & $p\,1\,q,r,s$ & $10.74025767\ldots
$ & $6V_1$
\\

$34)$ & $9_{28}^2$ & $2\,1,2,2++$ & $p\,1,q,r++$ & $10.74025767\ldots $ & $6V_1 $ \\
\hline

$35)$ & $9_2^3$ & $2\,1\,1\,1,2,2$ & $p\,1\,1\,1,q,r$ &
$11.76223429\ldots $ & $19.5826692\ldots $
\\

$36)$ & $9_{26}^2$ & $2\,1\,1,2,2+$ & $p\,1\,1,q,r+$ & $11.76223429\ldots $ & $19.5826692\ldots $ \\
\hline

$37)$ & $9_4^3$ & $2\,1,2,2,2$ & $p\,1,q,r,s$ & $12.2765628\ldots $
& $24.55255516\ldots $
\\

$38)$ & $9_{30}^2$ & $(2\,1,2)\,(2,2)$ & $(p\,1,q)\,(r,s)$ & $12.2765628\ldots $ & $24.55255516\ldots $ \\
\hline

$39)$ & $9_{30}$ & $2\,1\,1,2\,1,2$ & $p\,1\,1\,q,\,1,r$ & $11.95452697\ldots $ & $19.58266925\ldots $ \\
\hline

$40)$ & $9_1^4$ & $2,2,2,2+$ & $p,q,r,s+$ & $11.75183617\ldots $ & $24.55255552\ldots $ \\

$41)$ & $9_8^3$ & $(2,2+)\,(2,2)$ & $(p,q+)\,(r,s)$ &
$11.75183617\ldots $ & $24.55255552\ldots $
\\

$42)$ & $9_9^3$ & $(2,2)\,1\,(2,2)$ & $(p,q)\,1\,(r,s)$ & $11.75183617\ldots $ & $24.55255552\ldots $ \\
\hline

$43)$ & $9_{25}^2$ & $2\,2,2,2+$ & $p\,q,r,s+$ & $11.38178609\ldots $ & $23.32819345\ldots $ \\
\hline

$44)$ & $9_{28}$ & $2\,1,2\,1,2+$ & $p\,1,q\,1r+$ & $11.56317702\ldots $ & $18.83168337\ldots $ \\
\hline

$45)$ & $9_{10}^3$ & $.2\,1\,1$ & $.p\,1\,1$ & $13.32336092\ldots $ & $15.4156985\ldots $ \\
\hline

$46)$ & $9_{11}^3$ & $.2\,1:2$ & $.p\,1:q$ & $13.04040137\ldots $ &
$17.47714082\ldots $
\\

$47)$ & $9_{38}^2$ & $.2\,1:2\,0$ & $.p\,1:q\,0$ & $13.04040137\ldots $ & $17.47714082\ldots $ \\
\hline

$48)$ & $9_{33}$ & $.2\,1.2$ & $.p\,1.q$ & $13.28045564\ldots $ & $18.10505153\ldots $ \\
\hline

$49)$ & $9_{32}$ & $.2\,1.2\,0$ & $.p\,1.q\,0$ & $13.09989985\ldots $ & $18.1050515\ldots $ \\
\hline

$50)$ & $9_{29}$ & $.2.2\,0.2$ & $.p.q\,0.r$ & $12.20585617\ldots $ & $19.3538168\ldots $ \\
\hline

$51)$ & $9_{41}^2$ & $2:2\,0:2\,0$ & $p:q\,0:r\,0$ & $12.95742943\ldots $ & $22.07666239\ldots $ \\
\hline

$52)$ & $9_{41}$ & $2\,0:2\,0:2\,0$ & $p\,0:q\,0:r\,0$ & $12.09893603\ldots $ & $21.1717152\ldots $ \\
\hline

$53)$ & $9_{38}$ & $.2.2.2$ & $.p.q.r$ & $12.9328587\ldots $ & $20.72523729\ldots $ \\
\hline

$54)$ & $9_{40}^2$ & $2:2:2$ & $p:q:r$ & $12.04609204\ldots $ & $18.8316834\ldots $ \\
\hline

$55)$ & $9_{39}^2$ & $.2.2.2\,0$ & $.p.q.r\,0$ & $12.53617026\ldots $ & $19.7968546\ldots $ \\
\hline

$56)$ & $9_{39}$ & $2:2:2\,0$ & $p:q:r\,0$ & $12.81031\ldots $ & $21.0293868\ldots $ \\
\hline

$57)$ & $9_{12}^3$ & $.(2,2)$ & $.(p,q)$ & $13.81327844\ldots $ & $19.66433108\ldots $ \\
\hline

$58)$ & $9_{42}^2$ & $8^*2$ & $8^*p$ & $13.9484177\ldots $ & $16.0562293\ldots $ \\
\hline

$59)$ & $9_{34}$ & $8^*2\,0$ & $8^*p\,0$ & $14.34458139\ldots $ & $16.69568447\ldots $ \\
\hline

\end{tabular}

\bigskip

\begin{tabular}{|c|c|c|c|c|c|} \hline

$60)$ & $6_3^3$ & $2,2,-2$ & $p,q,-r$ & $0$ & $4V_1$ \\
\hline

$61)$ & $7_8^2$ & $2\,1,2,-2$ & $p\,1,q,-r$ & $V_1$ & $16.0004687\ldots $ \\
\hline

$62)$ & $8_{21}$ & $2\,1,2\,1,-2$ & $p\,1,q\,1,-r$ & $6.78371352\ldots $ & $17.6277542\ldots $ \\
\hline

$63)$ & $8_{15}^2$ & $2\,2,2,-2$ & $p\,q,r,-s$ & $V_1$ & $6V_1$ \\
\hline

$64)$ & $8_{16}^2$ & $2\,1\,1,2,-2$ & $p\,1\,1,q,-r$ & $5.3334895669\ldots $ & $17.6277542\ldots $ \\
\hline

$65)$ & $8_2^4$ & $2,2,2,-2$ & $p,q,r,-s$ & $2V_1$ & $24.09218408\ldots $ \\

$66)$ & $8_9^3$ & $(2,2)\,(2,-2)$ & $(p,q)\,(r,-s)$ & $2V_1$ & $24.09218408\ldots $ \\
\hline

$67)$ & $8_3^4$ & $2,2,-2,-2$ & $p,q,-r,-s$ & $0$ & $24.09218408\ldots $ \\

$68)$ & $8_{10}^3$ & $(2,2)\,-(2,2)$ & $(p,q)\,-(r,s)$ & $0$ & $24.09218408\ldots $ \\
\hline

$69)$ & $9_{44}$ & $2\,2,2\,1,-2$ & $p\,q,r\,1,-s$ & $7.4067675724\ldots $ & $23.32819345\ldots $ \\
\hline

$70)$ & $9_{45}$ & $2\,1\,1,2\,1,-2$ & $p\,1\,1,q\,1,-r$ & $8.6020031166\ldots $ & $19.58266925\ldots $ \\
\hline

$71)$ & $9_{48}^2$ & $2\,2\,1,2,-2$ & $p,q\,1,r,-s$ & $7.706911803\ldots $ & $19.5826692\ldots $ \\

$72)$ & $9_{14}^3$ & $2\,1\,1\,1,2,-2$ & $p\,1\,1\,1,q,-r$ & $7.706911803\ldots $ & $19.5826692\ldots $ \\
\hline

$73)$ & $9_{13}^3$ & $2\,1\,2,2,-2$ & $p\,1\,q,r,-s$ & $5.3334895669\ldots $ & $6V_1$ \\
\hline

$74)$ & $9_{16}^3$ & $2\,1,2,2,-2$ & $p\,1,q,r,-s$ & $9.966511884\ldots $ & $24.55255516\ldots $ \\

$75)$  & $9_{58}^2$ & $(2\,1,-2)\,(2,2)$ & $(p\,1,-q)\,(r,s)$ & $9.966511884\ldots $ & $24.55255516\ldots $ \\
\hline

$76)$ & $9_{56}^2$ & $(2\,1,2)\,(2,-2)$ & $(p\,1,q)\,(r,-s)$ & $8.997351944\ldots $ & $24.55255516\ldots $ \\
\hline

$77)$  & $9_{60}^2$ & $(2\,1,2)\,-(2,2)$ & $(p\,1,q)\,-(r,s)$ & $5.333489567\ldots $ & $24.55255516\ldots $ \\
\hline

$78)$ & $9_{18}^3$ & $(2,2+)\,(2,-2)$ & $(p,q+)\,(r,-s)$ & $2V_1$ & $24.55255516\ldots $ \\

$79)$ & $9_{19}^3$ & $(2,2+)\,-(2,2)$ & $(p,q+)\,-(r,s)$ & $2V_1$ & $24.55255516\ldots $ \\
\hline

$80)$ & $9_{61}^2$ & $2:-2\,0:-2\,0$ & $p:-q\,0:-r\,0$ & $0$ & $22.07666239\ldots $ \\
\hline

$81)$ & $9_{49}$ & $-2\,0:-2\,0:-2\,0$ & $-p\,0:-q\,0:-r\,0$ & $9.427073628\ldots $ & $21.1717152\ldots $ \\
\hline

$82)$ & $9_{20}^3$ & $.(2,-2)$ & $.(p,-q)$ & $3V_1$ & $19.66433108\ldots $ \\
\hline

$83)$ & $9_{21}^3$ & $.-(2,2)$ & $.-(p,q)$ & $0$ & $21.1717152\ldots $ \\
\hline

$84)$ & $9_{47}$ & $8^*-2\,0$ & $8^*-p\,0$ & $10.0499579\ldots $ & $16.69568447\ldots $ \\
\hline
\end{tabular}

\bigskip

\normalsize

Some of non-alternating source links from the preceding table are
non-hyperbolic, so their hyperbolic volume is 0. In the family 60),
link $p,q,-r$ is non-hyperbolic for $r=2$, and hyperbolic otherwise;
in the family 66) all links except source link $2,2,-2,-2$ are
hyperbolic. The same holds for the family 67) and its source link
$(2,2)\,-(2,2)$, and for family 80) and its source link
$2:-2\,0:-2\,0$.

It is interesting to notice that all families with the same lower
and upper bound of hyperbolic volume have subfamilies of distinct
links with the same hyperbolic volume, which can be distinguished by
Alexander and Jones polynomial. Alternating link families 3) and 4)
have the subfamilies $2\,p\,2$ and $p,2,2$ of links which cannot be
distinguished by hyperbolic volume, 8) and 9) have the subfamilies
$p\,1,2,2$ and $p,2,2+$, 14) and 15) the subfamilies $p,q,2,2$ and
$(p,q)\,(2,2)$, 18) and 19) the subfamilies $p\,1\,1,2,2$ and
$p\,1,2,2+$, 22) and 23) the subfamilies $.p:2$ and $.p:2\,0$, 31)
and 32) the subfamilies $p\,1\,q,2,2$ and $p\,1,q,2++$, 35) and 36)
the subfamilies $p\,1\,1\,1,2,2$ and $p\,1\,1,2,2+$, 35) and 36) the
subfamilies $p\,1,q,2,2$ and $(p\,1,q)\,(2,2)$, 40), 41) and 42) the
subfamilies $(2,2)\,p\,(2,2)$, $2,2,2,2+p$ and
$(2,2+p)\,(2,2)$\footnote{Sequence of $p$ pluses is denoted by
$+p$.}, and 46) and 47) have the subfamilies $.p\,1:2$ and
$.p\,1:2\,0$ with the same property. The same holds for
non-alternating link families 65) and 66) and their subfamilies
$p,q,2,-2$ and $(p,q)\,(2,-2)$, 71) and 72) and their subfamilies
$p\,2\,1,2,-2$ and $p\,1\,1\,1,2,-2$, 74) and 75) and their
subfamilies $p\,1,2,2,-2$ and $(p\,1,-2)\,(2,2)$, and 78) and 79)
and their subfamilies $(2,2+)\,(p,-2)$ and $(2,2+)\,-(p,2)$. Even in
the limiting case, for the families 67) and 68) with non-hyperbolic
source links and with the same upper bound, their subfamilies
$p,q,-2,-2$ and $(p,q)\,-(2,2)$ consist of links with the same
hyperbolic volume. We propose the following conjecture:

\medskip

\noindent {\bf Conjecture 2.} Every two families with the same lower
and upper bound of hyperbolic volume have subfamilies of links with
the same hyperbolic volumes.

\medskip

\begin{theorem}
Hyperbolic volume completely distinguishes alternating links
belonging to the same family.
\end{theorem}

\medskip

This theorem does not hold for non-alternating links: for example,
links of the subfamilies $p,3,-2$ and $p-6,2,-3$ ($p\ge 8$)
belonging to the family 60) $p,q,-r$ have the same hyperbolic
volume. These links can be distinguished by Alexander and Jones
polynomial.

For $n\le 12$ there are no source knots with the same hyperbolic
volume, so we propose the following conjecture:

\bigskip

\noindent {\bf Conjecture 3.} Hyperbolic volume completely
distinguishes source knots\footnote{Except mutant source knots.}.

\medskip

\begin{definition}
The replacement of a tangle $(2,2)$ by $(2,2)\,0$, or $(2,-2)$ by
$(2,-2)\,0$, or $-(2,2)$ by $-(2,2)\,0$ and {\it vice versa} will be
called {\it $(2,2)$-reversal}. Two links are called
$(2,2)$-equivalent if one can be obtained from the other by
$(2,2)$-reversals.
\end{definition}

\medskip

\begin{theorem}
$(2,2)$-equivalent links have the same hyperbolic volume. The same
holds for their corresponding augmented links.
\end{theorem}

\medskip

From the Theorem 3 we can make many conclusions about links and
their hyperbolic volume. For example, pretzel links $2,2,-p$ ($p\ge
3$) and rational links $2\,(p-2)\,2$ have the same hyperbolic volume
and its upper bound $4V_1$. We can also conclude that all links
belonging to the family $2,2,-p$ are hyperbolic, except the link
$2,2,-2$ which is $(2,2)$-equivalent with non-hyperbolic rational
link $2\,-2\,2=4$.

\bigskip
\bigskip

{\bf References}

\medskip

\noindent [1] Jablan, S.~V. and Sazdanovi\' c, R. (2007) {\it
LinKnot- Knot Theory by Computer}. World Scientific, New Jersey,
London, Singapore.

\medskip

\noindent [2] Stoimenow, A.: Graphs, determinants of knots and
hyperbolic volume.

\medskip

\noindent [3] Dasbach, O.~T. and X-S. Lin (2004) A volume-ish
theorem for the Jones polynomial of alternating knots,
arXiv:math/0403448v1 [math.GT]

\medskip

\noindent [4] Petronio,~C. and Vesnin,~A. (2007) Two-sided bounds
for the complexity of cyclic branched coverings of two-bridge links,
arXiv:math/0612830v2 [math.GT]

\medskip

\noindent [5] Lackenby,~M., Agol,~I and Thurston,~D. (2004) The
volume of hyperbolic alternating link complements, Proc. London
Math. Soc. (3) {\bf 88} 204--224.

\medskip

\noindent [6] Purcell,~J. (2007) Volumes of highly twisted knots and
links, arXiv:math/060447v2 [math.GT]

\medskip

\noindent [7] Futer,~D., Kalfagianni,~E. and Purcell, J. (2008) Dehn
filling, volume, and the Jones polynomial, arXiv:math/0612138v4
[math.GT]

\medskip

\noindent [8] Rolfsen, D. (1976) Knots and Links, Publish \& Perish
Inc., Berkeley (American Mathematical Society, AMS Chelsea
Publishing, 2003).

\end{document}